\documentclass[12pt,a4paper]{amsart}
\usepackage{amsfonts,amssymb}
\usepackage[english]{babel}
\pdfoutput=1
\title
{On orthogonal systems in Hilbert $C^*$-modules}
\date{4 June 2009}

\author{Giovanni Landi}
\address{Dipartimento di Matematica e Informatica, Universit{\`a} di Trieste}
\email{landi@univ.trieste.it}
 
\author{Alexander Pavlov}
\address{Dipartimento di Matematica e Informatica, Universit{\`a} di
Trieste} \email{axpavlov@gmail.com}

\newtheorem{teo}{Theorem}[section]
\newtheorem{lem}[teo]{Lemma}
\newtheorem{prop}[teo]{Proposition}

\newtheorem{cor}[teo]{Corollary}

\newtheorem{ex}[teo]{Example}

\begin{document}

\begin{abstract}  Analogues for Hilbert $C^*$-modules of classical results of Fourier series
theory in Hilbert spaces are considered. Relations between
different properties of orthogonal and orthonormal systems for
Hilbert $C^*$-modules are studied with special attention paid on
the differences with the well-known Hilbert
space situation. \\
\end{abstract}

 \maketitle


\section{Introduction}
In this paper we study  properties of orthogonal and orthonormal
systems in Hilbert $C^*$-modules. Actually the theory of Hilbert
$C^*$-modules is at an intermediate stage between the theory of Hilbert spaces
and the theory of general Banach spaces and can be considered as a
`quantization' of the Hilbert space theory. Roughly speaking by
quantization here we mean the following: there are crucial notions
and definitions of the theory that include commutative
objects like functions or just scalars and one replaces them in some
proper way by noncommutative objects like elements of an
arbitrary $C^*$-algebra. From this point of view the definition of
Hilbert $C^*$-modules can be obtained by replacing complex
vector spaces with modules over a $C^*$-algebra and allowing the
inner product to take values in this $C^*$-algebra. This concept
originally arose in \cite{Kap52} for commutative $C^*$-algebras
and it was studied in the general noncommutative context in
\cite{Pas73, Rieffel74}. The theory of Hilbert $C^*$-modules has a
number of effects, related to the operator nature of the
`coefficients' of their elements, that make it much more complicated to handle
with respect to
the usual Hilbert space theory. For example, a closed submodule of a
Hilbert $C^*$-module need not be orthogonally (or even
topologically - in the sense of direct sums of closed Banach
submodules) complemented, a bounded $A$-linear operator in a
Hilbert module over a $C^*$-algebra $A$ need not have and adjoint, a Hilbert $C^*$-module
need not be self-dual, i.e. canonically isomorphic to its $C^*$-dual module
(cf. \cite{Pas73, Lance, MaTrobook}).

 Any Hilbert space can be described as a space of sequences
(or nets in the non-separable case) $\{c_i\}$ of complex numbers
such that the series $\sum\nolimits_i c_i^*c_i$ converge in norm.
The reason that any vector is represented in a unique way by its
coordinate sequence is explained via the Fourier series theory:
any Hilbert space admits a complete orthonormal system which
automatically has to be closed (this exactly means that the
Parseval equality is valid for the system); consequently, it forms
an orthonormal basis for the Hilbert space. Unfortunately -- but
not surprisingly --, for Hilbert $C^*$-modules this scheme does
not work and they do not admit orthonormal bases in general (e.g.
\cite{Lance, MaTrobook}). The reason is, as we will discuss
more thoroughly below, that the Fourier series of a vector
$x$ of a Hilbert $C^*$-module $M$ with respect to a certain
orthonormal system of $M$ need not converge in norm to $x$ even
when this orthonormal system is complete (Example
\ref{ex:ONS_full_but_not_closed}).

An efficient way to cope with this difficulty is provided by the
concept of frame that was introduced in \cite{FrankLar,
FrankLarJOP} for countably generated Hilbert modules. We remind
that a sequence $\{x_i\}$ of vectors of a Hilbert module over a
unital $C^*$-algebra is called a \emph{frame} if for any vector
$x$ in the Hilbert module there are real constants $C, D>0$ such
that
\[
C\langle x,x\rangle\le \sum\nolimits_i \langle x,x_i\rangle\langle
x_i,x\rangle\le D\langle x,x\rangle.
\]
The frame is said to be \emph{tight} if $C=D$, and it is said to
be \emph{normalized} if $C=D=1$.  The frame is named
\emph{standard} if one has that $\langle x,x\rangle =
\sum\nolimits_i \langle x,x_i\rangle\langle x_i,x\rangle $ for any
vector $x$ in the Hilbert module, a condition which is the
analogue of the Parseval equality. The following crucial result
about frames describes the conditions so that the reconstruction
formula holds.

\begin{teo}\label{teo:FLJOP_th4.1}
 {\rm (cf. \cite[Th. 4.1]{FrankLarJOP})} Let $A$ be a unital
$C^*$-algebra, $M$ be a finitely or countably generated Hilbert
$A$-module and $\{x_i\}$ be a normalized tight frame of $M$. Then
the reconstruction formula
\[
x=\sum\nolimits_{i} e_i \langle e_i,x\rangle
\]
holds for every $x\in M$, in the sense of convergence in norm,  if
and only if the frame $\{x_i\}$ is standard.
\end{teo}

Just to mention some applications, the  frame  approach has
already shown its usefulness for the description of conditional
expectations of finite index and for the analysis of some classes
of $C^*$-algebras (see references in \cite{FrankLar,
FrankLarJOP}). It is also very useful to investigate finitely
generated projective modules arising from vector bundles and in
particular for finding bases for the space of sections of
non-trivial vector bundles \cite[Proposition 7.2]{Rieffel06}.

In the present paper we will not deal directly with frames, but
nevertheless our considerations are very close to the frame
approach. Also, as it will be discussed thoroughly in
\S\ref{se:obshm}, some help to overcome the lack of orthonormal
bases for a general Hilbert module sometimes comes from the
Kasparov's stabilization theorem \cite{KaspJO}. Our aim with the
present note is two-fold. On the one hand we seek to obtain
natural analogues -- for arbitrary (i.e. non necessarily countably
generated) Hilbert modules over operator algebras -- of well-known
results about Fourier series and orthonormal systems in Hilbert
spaces. On the other hand to highlight, mainly using examples,
some of the differences between these two theories. We will show
that any Hilbert $C^*$-modules have complete orthogonal systems
(Proposition \ref{teo:full_orth_syst}), but a complete orthogonal
and even orthonormal system needs not be closed at the same time
(Examples \ref{ex:x_in_C_0[0,1]},
\ref{ex:ONS_full_but_not_closed}). Also the completeness  of an
orthogonal system does not imply that it forms a basis even in
some weak sense (Example \ref{ex:full_without_(3)}); despite these
results, there is an analogue of the Bessel inequality for Hilbert
$C^*$-modules. Fourier series of vectors with respect to some
orthogonal system need not converge in norm, but only with respect
to the strong topology (Theorem \ref{teo:strong_convergence_FS}).
We also describe interrelations between different properties of orthonormal
systems in Hilbert $C^*$-modules (Theorem \ref{teo:equiv_cond_for_ONS}, Corollary
\ref{cor:equiv_cond_for_ONS}).

\section{Orthogonal systems in Hilbert $C^*$-modules}
In the sequel, $(M,\langle\cdot,\cdot\rangle)$ is always a Hilbert
module over a $C^*$-algebra $A$, unless otherwise
explicitly stated.  A collection $\{e_i\}_{i\in I}$, indexed by some set
$I$, of vectors from $M$ is called \emph{orthogonal} if $\langle
e_i,e_j\rangle=0$ whenever $i\ne j$. The orthogonal system
$\{e_i\}_{i\in I}$ is said to be \emph{quasi-orthonormal} if there
are (self-adjoint) projections $p_i$ in $A$ such that $\langle
e_i,e_i\rangle=p_i$ for all $i\in I$ and it is said to be
\emph{orthonormal} provided $A$ is unital (where this not the case
we would be able to join the unit, but there will be no need for
such complications in the following) and for the inner squares it happens that
$\langle e_i,e_i\rangle=1$ for all $i\in I$.

Let $\{e_i\}_{i\in I}$ be an orthogonal system of $M$, $x$ be an
arbitrary vector in $M$ and $F\subset I$ be any finite subset. Then
\begin{gather*}
    S_F=\sum\nolimits_{i\in F} e_i \langle e_i,x\rangle
\end{gather*}
stands for the corresponding partial sum of the Fourier series
with respect to $\{e_i\}_{i\in I}$ and a straightforward
computation provides the formula:
\begin{gather}\label{eq:analogue_Pars_eq}
\langle x-S_F, x-S_F\rangle=\langle x,x\rangle-
     2\sum\nolimits_{i\in F}\langle x,e_i\rangle\langle e_i,x\rangle+
      \sum\nolimits_{i \in F}\langle x, e_i\rangle \langle e_i, e_i\rangle
                    \langle e_i, x\rangle.
\end{gather}

\medskip
Next, given an orthogonal system $\{e_i\}_{i\in I}$ of $M$, in the
sequel we will explore consequences and relations among them of the following conditions:
\begin{itemize}
    \item[{\bf (c1)}] \label{eq:basis1}
    The system $\{e_i\}_{i\in I}$ generates $M$ over $A$,
    \begin{gather*}
    M=\overline{\mathrm{span}_A\{e_i : i\in I\}} ,
    \end{gather*}
 that is to say, the closure of its $A$-linear span coincides with $M$.
       \item[]
    \item[{\bf (c2)}] \label{eq:basis2}
    For any $x$ of $M$ there are elements $a_i$ of $A$ such that
    \begin{gather*}
     x=\sum\nolimits_{i\in I} e_i a_i,
    \end{gather*}
    where convergence in norm is meant and
    \[
\sum\nolimits_{i\in I} e_i a_i=\lim_{F\in \mathcal{F}}\sum\nolimits_{i\in F} e_i a_i
    \]
    indicates the limit over the set $\mathcal{F}$ of all finite
    subsets of $I$, directed by inclusions.
           \item[]
    \item[{\bf (c3)}] \label{eq:basis5}
The system $\{e_i\}_{i\in I}$ is said to be \emph{closed} if it
happens that for any $x\in M$ the series
\begin{gather*}
\sum\nolimits_{i\in I}\bigl( 2\langle x,e_i\rangle\langle e_i,x\rangle-
 \langle x,e_i\rangle\langle e_i,e_i\rangle\langle e_i,x\rangle \bigr )
\end{gather*}
converges in norm to $\langle x,x\rangle$. Using
(\ref{eq:analogue_Pars_eq}) this  exactly means that any vector of
$M$ is the limit in norm of its Fourier series.
\item[]
    \item[{\bf (c4)}] \label{eq:basis6}
    The system $\{e_i\}_{i\in I}$ is said to be \emph{complete}
provided there is no non-zero vector $x$ of $M$ such that $\langle
e_i, x\rangle=0$ for all $i\in I$.
    \end{itemize}

\bigskip
It is clear  that condition {\bf (c2)} implies {\bf (c1)}.
The next example shows that the converse is not true: condition
{\bf (c1)} does not imply {\bf (c2)} in general.

\begin{ex}\label{ex:x_in_C_0[0,1]}\rm
Let $A={ C_0(0,1]}=\{f\in C[0,1] : f(0)=0\}$ and $M=A$ be the
Hilbert $A$-module with respect to the inner product:
\begin{gather*}
\langle a,b\rangle=a^*b,\quad a,b\in A.
\end{gather*}
Then, the one-point-set $\varepsilon=\{f\}$, where $f(x)=x$ for
$x\in [0,1]$, is an orthogonal and, clearly, complete system of
$M$. Suppose $B$ is the closure of the $*$-algebra
\begin{gather*}
    \{fg : g\in {C_0(0,1]}\}.
\end{gather*}
Then, with $\{g_i\}$ standing for the approximative identity of ${
C_0(0,1]}$, the $C^*$-algebra $B$ contains $f$ as the limit
$f=\lim_i fg_i$. As a consequence, $B$ separates points of the
interval $[0,1]$ and, consequently, coincides with ${ C_0(0,1]}$
by the Stone-Weierstrass theorem (cf. \cite[Theorem
IV.10]{Reed_Simon}). Thus, the system $\varepsilon$ satisfies {\bf
(c1)}. But at the same time it does not satisfy {\bf (c2)} since,
for instance, the function $f$ cannot be represented as a product
$fg$ for any $g\in { C_0(0,1]}$.
\end{ex}

The next example shows there are complete orthogonal systems not
satisfying {\bf (c1)}.

\begin{ex}\label{ex:full_without_(3)}\rm We will slightly modify Example
\ref{ex:x_in_C_0[0,1]}. Let $A={ C_0(0,1]}$ and $M=A$ again, but
now take $\varepsilon=\{g\}$, where
\[
g(x)=
\begin{cases}
x, & \text{if $\,\, 0\le x\le 1/2$;}\\
1-x, & \text{if $\,\, 1/2< x\le 1$.}
\end{cases}
\]
Clearly, $\varepsilon$ is a complete orthogonal system for $M$.
But it cannot satisfy {\bf (c1)} since the closure of the set
$\{gh : h\in { C_0(0,1]}\}$ belongs to the suspension $SA=\{f\in
C[0,1] : f(0)=f(1)=0\}$ of $A$ rather then to $A$ itself.
\end{ex}

A use of  the Zorn lemma directly
ensures that any pre-Hilbert $C^*$-module admits a complete
orthogonal system, more precisely the next statement is true.

\begin{prop}\label{teo:full_orth_syst}
Every orthogonal system of a pre-Hilbert $C^*$-module can be
enlarged to a complete orthogonal system; also, every orthogonal
system of norm one vectors can be enlarged to a complete
orthogonal system of norm one vectors. $\;\Box$
 \end{prop}

This observation may be strengthen for von Neumann modules (see
\cite{SkeideRIA}). Indeed, let $B(G)$ denotes the set of all
linear bounded operators in a Hilbert space $G$, $A\subset B(G)$
be a von Neumann algebra acting non-degenerately on G,  and $M$ be
a Hilbert $A$-module. Then the algebraic tensor product $M\otimes
G$ becomes a pre-Hilbert space with respect to the inner product
$\langle x\otimes g,x'\otimes g'\rangle=\langle g,\langle
x,x'\rangle g'\rangle$. Let $H=\overline{M\otimes G}$ stands for
the Hilbert space completion of $M\otimes G$. We can consider in a
natural way the module $M$ as a linear subspace of the space
$B(G,H)$ of all bounded linear operators from $G$ to $H$. Then $M$
is said to be a \emph{von Neumann module} if it is strongly closed
in $B(G,H)$. These modules behave themselves like Hilbert spaces,
mostly because they are necessarily self-dual. As for the fact
that any von Neumann module admits a complete quasi-orthonormal
system one has the following (\cite[Theorem 4.11]{SkeideRIA}).

Now, the analogue of the Bessel inequality for an
orthogonal system $\{e_i\}_{i\in I}$ of norm one vectors in a
pre-Hilbert $C^*$-module exists only as a finite version,
i.e.
$$ \sum\nolimits_{i\in F} \langle x,e_i\rangle\langle e_i,x\rangle\le
                \langle x,x\rangle
                $$
holds for every finite subset $F\subset I$ and any vector $x$. But
provided $\{e_i\}_{i\in I}$ is made of norm one vectors and fulfils {\bf (c3)} the
restriction on finiteness may be omitted, i.e. under these additional
conditions
$$
\sum\nolimits_{i\in I} \langle x,e_i\rangle\langle
e_i,x\rangle\le \langle x,x\rangle.
$$
To make certain of the last inequality we need just  a direct use
of the following auxiliary result.
\begin{lem} Let $\{e_i\}_{i\in I}$ be an orthogonal system of
norm one vectors in a pre-Hilbert $C^*$-module. Then the following
conditions are equivalent:
\begin{enumerate}
    \item the net
    $$
    \left\{A_F=\sum\nolimits_{i\in F}\left(2\langle x,e_i\rangle\langle e_i,x\rangle-
 \langle x,e_i\rangle\langle e_i,e_i\rangle\langle e_i,x\rangle\right) : F \,
 \text{is a finite subset of}\,\, I\right\}
 $$
 converges in norm;
    \item the net $\left\{B_F=\sum\nolimits_{i\in F}\langle x,e_i\rangle\langle e_i,x\rangle
  : F \, \text{is a finite subset of}\,\, I\right\}$ converges in
 norm.
\end{enumerate}
\end{lem}

\begin{proof} Clearly, \rm{(ii)} implies \rm{(i)}, so we just need
 verify the inverse implication. Assume \rm{(i)} is true and
 denote $C_F=\sum\nolimits_{i\in F} \langle x,e_i\rangle\langle e_i,e_i\rangle\langle
 e_i,x\rangle$. Then $B_F=A_F-B_F+C_F$, whence $B_F$ converges if and
 only if $B_F-C_F$ does. To finish the argument it only remains
 to observe that  $B_F-C_F\le A_F$.
\end{proof}

Although the Parseval equality does not take holds for arbitrary
Hilbert $C^*$-modules, there is a weakened version in the
$W^*$-case.

\begin{prop}\label{teo:strong_convergence_FS}
Suppose $\{e_i\}_{i\in I}$ is an orthogonal system of norm one
vectors in a Hilbert module $M$ over a von Neumann algebra. Then for
any vector $x\in M$ the net
\[
a_F=2\sum\nolimits_{i\in F} \langle x, e_i\rangle\langle e_i, x\rangle-
\sum\nolimits_{i\in F} \langle x, e_i\rangle \langle e_i, e_i\rangle
\langle e_i, x\rangle,
\]
indexed by finite subsets $F$ of $I$, converges with respect
to the strong topology.
\end{prop}

\begin{proof}
Clearly, the elements $a_F$ are positive, and the equality
(\ref{eq:analogue_Pars_eq}) implies $a_F\le \langle x,x\rangle$
for all finite subsets $F$ of $I$. It only remains to check that
the net $\{a_F\}$ is not decreasing, and the required result
will follow from \cite[Theorem 4.1.1]{Murphy}. So, let $G$, $F$ be
finite subsets of $I$ and $F\subset G$;  one gets:
\begin{align*}
a_G-a_F&=2\sum\nolimits_{i\in G\setminus F}\langle x,e_i\rangle\langle
e_i,x\rangle- \sum\nolimits_{i\in G\setminus F} \langle x, e_i\rangle
\langle e_i, e_i\rangle \langle e_i, x\rangle\\
&=\sum\nolimits_{i\in G\setminus F}\langle x,e_i\rangle\langle
e_i,x\rangle+ \sum\nolimits_{i\in G\setminus F} \langle e_i, x\rangle^*
(1-\langle e_i, e_i\rangle) \langle e_i, x\rangle\ge 0 ,
\end{align*}
under our assumption that the $e_i$'s are norm one vectors.
\end{proof}

\begin{lem}\label{teo:optimality_property}
\rm{(The optimality property of Fourier series)}. Suppose
$\{e_i\}_{i\in I}$ is an orthonormal system of $M$, $x\in M$ is an
arbitrary vector. Then
\begin{gather*}
\langle x-S_F, x-S_F\rangle\le
     \Bigl\langle x-\sum\nolimits_{i\in F}e_ia_i \,, x-\sum\nolimits_{i\in F}e_ia_i\Bigr\rangle
\end{gather*}
for any elements $a_i\in A$ and for any finite subset $F\subset
I$. Moreover, in the above expression the equality occurs if and only if
$a_i=\langle e_i,x\rangle$ for any $i\in F$.  
\end{lem}

\begin{proof}
The above inequality follows from the following sequence of
transformations:
\begin{multline*}
\langle x-\sum\nolimits_{i\in F} e_i a_i \,, x-\sum\nolimits_{i\in F} e_i a_i\rangle  =
  \langle x,x\rangle-
  \sum\nolimits_{i\in F} \langle x,e_i\rangle a_i-
  \sum\nolimits_{i\in F} a_i^*\langle e_i,x\rangle  +
  \sum\nolimits_{i\in F} a_i^*a_i \\
  =
\langle x,x\rangle-
  \sum\nolimits_{i\in F} \langle x, e_i\rangle\langle e_i, x\rangle+
  \sum\nolimits_{i\in F}\left(a_i-\langle e_i,x\rangle\right)^*
               \left(a_i-\langle e_i,x\rangle\right) \\
   =
 \langle x-S_F, x-S_F\rangle+
  \sum\nolimits_{i\in F}\left(a_i-\langle e_i,x\rangle\right)^*
               \left(a_i-\langle e_i,x\rangle\right),
\end{multline*}
using identity (\ref{eq:analogue_Pars_eq}) for orthonormal systems: $\langle x-S_F,
x-S_F\rangle=\langle x, x\rangle-
     \sum\nolimits_{i\in F} \langle x, e_i\rangle\langle e_i, x\rangle$.
\end{proof}

The next example emphasizes that for bases which are orthogonal
but not orthonormal there is no uniqueness of the decomposition
{\bf (c2)}.

\begin{ex}\rm \label{ex:orthogonal_basis}
Consider $A=L^\infty [0,1]$, $M=A$ with the usual inner product and let an orthogonal system of $M$ be given by $\varepsilon=\{f_1,f_2\}$, where
\[
f_1(x)=
\begin{cases}
1, & \text{if $\,\, 0\le x\le 1/2$;}\\
0, & \text{if $\,\, 1/2< x\le 1$}
\end{cases}
\]
and
\[
f_2(x)=
\begin{cases}
0, & \text{if $\,\, 0\le x\le 1/2$;}\\
1, & \text{if $\,\, 1/2< x\le 1$.}
\end{cases}
\]
Then $g=f_1g+f_2g$ for any $g\in A$, so $\varepsilon$ forms the
basis in the sense of the condition {\bf (c2)}. But now
uniqueness does not hold since, for instance, the unit function
$\mathbf{1}$ of $A$ complies:
\[
\mathbf{1}=f_1\cdot\mathbf{1}+f_2\cdot\mathbf{1}=f_1\cdot
f_1+f_2\cdot f_2.
\]
\end{ex}

\bigskip
     Using a Banach space like terminology (cf.
\cite{LT_CBS}) we say that an orthogonal system $\{e_i\}_{i\in I}$
of $M$ forms an orthogonal \emph{Schauder basis} for $M$ (over
$A$) if $\{e_i\}_{i\in I}$ satisfies {\bf (c2)} and the
coefficients in the decomposition {\bf (c2)} are unique for any
vector $x$ of $M$.

Let us remind (cf.~\cite{MaTrobook}) that an element $x$ of $M$ is
called \emph{non-singular} if its inner square $\langle
x,x\rangle$ is invertible in $A$.  Clearly, an orthogonal system
$\{e_i\}_{i\in I}$ of  $M$ satisfying the condition {\rm {\bf
(c2)}} is an orthogonal Schauder basis provided it consists of
non-singular vectors. Indeed, in this case the coefficients $a_i$
of the decomposition {\bf (c2)} take the form
\[
a_i=\langle e_i,e_i\rangle^{-1}\langle e_i,x\rangle,
\]
from which one infers their uniqueness. The next theorem gives
additional properties.

\begin{teo}\label{teo:singular vectors}
 Assume an orthogonal system $\{e_i\}_{i\in I}$ of a
Hilbert module $M$ over an unital $C^*$-algebra $A$, satisfying
the condition {\rm {\bf (c2)}}, contains at least one singular
vector, $e_t say$. Then the system $\{e_i\}_{i\in I}$ does not
form a Schauder basis if at least one of the following conditions
holds:
\begin{enumerate}
    \item  zero is an isolated
    point of the spectrum of $\langle e_t,e_t\rangle$;
    \item for any element $a$ of
    $A$ which is both non-invertible and non-zero,
    there is a non-zero element $b$ of $A$ such that $ab=0$.
\end{enumerate}
\end{teo}

\begin{proof} Firstly,  suppose that (i) is true and consider the following
continuous function
\[
f(x)=
\begin{cases}
1, & \text{if $\,\, x=0$;}\\
0, & \text{otherwise} ,
\end{cases}
\]
on the spectrum of $\langle e_t,e_t\rangle$. Then, the element
$b=f(\langle e_t,e_t\rangle)$ is not zero
\cite[VII.3]{Reed_Simon}, belongs to $A$ and
\begin{gather}\label{eq:singular vectors}
\langle e_t,e_t\rangle \, b=0.
\end{gather}
Therefore $e_t (b+\mathbf{1})=e_t\mathbf{1}$ meaning that
$\{e_i\}_{i\in I}$ does not form a Schauder basis. The same
argument is valid under the assumption (ii) as well, because it
directly yields the equality~(\ref{eq:singular vectors}).
\end{proof}

\bigskip
We give examples of $C^*$-algebras with and without the property
(ii) of Theorem~\ref{teo:singular vectors}.

\begin{ex}\rm Any unital commutative $C^*$-algebra, for instance
$C[0,1]$, does not satisfy the condition (ii) in Theorem~\ref{teo:singular vectors}. 
Condition (ii)
holds for finitely dimensional $C^*$-algebras. But the algebra
$B(H)$ of bounded linear operators on a separable Hilbert space
$H$ does not enjoy (ii). To see this it suffices to take the
operator $a=\mathrm{diag}(1,\frac{1}{2},\dots,\frac{1}{n},\dots)$;
it is compact and not invertible, but there is no non-zero $b$ in
$B(H)$ such that $ab=0$.
\end{ex}

\begin{teo} \label{teo:equiv_cond_for_ONS}
Let $\{e_i\}_{i\in I}$ be an orthonormal system in a Hilbert
module $M$ over a unital $C^*$-algebra $A$. Then the conditions
${\rm {\bf (c1)}} - {\rm {\bf (c3)}}$ are equivalent and each of
them is strictly stronger than {\rm {\bf (c4)}}.
\end{teo}

\begin{proof} It is clear that {\bf (c3)} implies the completeness  of
$\{e_i\}_{i\in I}$. On the other hand Example
\ref{ex:ONS_full_but_not_closed} below ensures that {\bf (c4)}
does not imply {\bf (c3)}.
To show that $\mathrm{{\bf (c1)}}$ implies ${{\bf (c2)}}$
let us consider an arbitrary vector $x\in M$. Then
for any $\delta >0$ one can find a finite subset $G\subset I$ and
elements $a_i\in A$ such that
\begin{gather*}
\biggl\|x-\sum\nolimits_{i\in G}e_ia_i\biggr\|<\delta.
\end{gather*}
Now for any finite set $F\subset I$  containing $G$  put
\[
b_i=
\begin{cases}
a_i, & \text{if $\,\, i\in G$;}\\
0,& \text{if $\,\, i\in F\setminus G$}.
\end{cases}
\]
Applying Lemma \ref{teo:optimality_property} we conclude that
\begin{gather*}
\delta >\biggl\|x-\sum\nolimits_{i\in F} e_ib_i\biggr\|\ge
\biggl\|x-\sum\nolimits_{i\in F}e_i\langle e_i,x\rangle\biggr\|.
\end{gather*}
This just means that $\lim_F\sum\nolimits_{i\in F}e_i\langle e_i,x\rangle=x$.

It is clear that the decomposition of $x$ in {\bf (c2)} is unique
for any $x\in M$, besides the coefficients $a_i=\langle
e_i,x\rangle$, so {\bf (c2)} implies {\bf (c3)}. Besides this,
obviously, {\bf (c3)} implies {\bf (c1)}. This finishes the proof.
\end{proof}

    According to \cite{FrankLarJOP} we  call an orthogonal system $\{e_i\}_{i\in
    I}$ of $M$ an orthogonal \emph{standard Riesz basis} if it is a standard
    frame satisfying {\bf (c1)} and endowed with the
additional property that $A$-linear combinations
$\sum\nolimits_{j\in S} e_ja_j$ with coefficients $a_j\in A$ and
$S\subset I$ are equal to zero if and only if $e_ja_j$ equals zero
for any $j\in S$.

\begin{cor}\label{cor:equiv_cond_for_ONS}
 Let $\{e_i\}_{i\in I}$ be an orthonormal system in a Hilbert
module $M$ over a unital $C^*$-algebra $A$. Then the following
conditions are equivalent:
\begin{enumerate}
    \item $\{e_i\}_{i\in I}$ is a Schauder basis;
    \item $\{e_i\}_{i\in I}$ is a standard Riesz basis;
    \item $\{e_i\}_{i\in I}$ satisfies any of the  conditions
    ${\rm {\bf (c1)}} - {\rm {\bf (c3)}}$.
\end{enumerate}
\end{cor}

\begin{proof} Since the decomposition of $x$ in {\bf (c2)} is unique
for any $x\in M$, {$\bf (c2)$} holds for $\{e_i\}_{i\in I}$ if and
only if $\{e_i\}_{i\in I}$ forms a Schauder basis. Moreover,
clearly, ${\rm (i)}\Rightarrow {\rm (ii)}\Rightarrow {\bf (c2)}$.
\end{proof}

\section{Orthonormal bases and standard Hilbert $C^*$-modules}\label{se:obshm}

The \emph{standard Hilbert module} over a $C^*$-algebra $A$, which
is denoted by $l_2(A)$ or $H_A$, consists of all sequences $(a_i)$
of elements of $A$ such that the series
$\sum\nolimits_{i=1}^\infty a_i^*a_i$ converges in norm. The inner
product of elements $x=(a_i)$ and $y=(b_i)$ of $l_2(A)$ is given
by $\langle x,y\rangle=\sum\nolimits_{i=1}^\infty a_i^*b_i.$
According to the Kasparov's stabilization theorem \cite{KaspJO}
every countably generated Hilbert $C^*$-module is a direct summand
of $l_2(A)$.
The notion of a standard Hilbert $C^*$-module can be naturally
generalized for any cardinality in the following way.  Let $I$ be
an arbitrary set and $(a_i)_{i\in I}$ be a collection of elements
from $A$ indexed by $I$. Given that the collection $\mathcal{F}$
of finite subsets of $I$ is partially ordered by inclusions we
form a net $\{\sum\nolimits_{i\in F} a_i^*a_i : F\in\mathcal{F}\}$
of finite sums. If this net converges in norm we will declare by
definition that the series $\sum\nolimits_{i\in I} a_i^*a_i$
converges in norm. Then the Hilbert $A$-module $H_{A,I}$ is made
of all collections $(a_i)_{i\in I}$ of elements of $A$ such that
the series $\sum\nolimits_{i\in I} a_i^*a_i$ converges in norm
with the inner product of elements $x=(a_i)$ and $y=(b_i)$ of
$H_{A,I}$ given by $\langle x,y\rangle=\sum\nolimits_{i\in I}
a_i^*b_i$.

 It is a well-known fact that a Hilbert module $M$ over a unital
$C^*$-algebra $A$ possesses an orthonormal system $\{e_i\}_{i\in
I}$ that satisfies the condition {\rm {\bf (c2)}} (such collection
of vectors is said to be an \emph{orthonormal basis}) if and only
if $M$ is isomorphic to the standard $A$-module $H_{A,I}$. Let us
recall just the  sketch of the proof of this assertion. The
orthonormal basis $\{e_i\}_{i\in I}$ of $M$ is closed by Theorem
\ref{teo:equiv_cond_for_ONS}, consequently both $x=\sum_{i\in
I}e_i\langle e_i,x\rangle$ for any $x$ of $M$ and the series
$\sum_{i\in I} \langle e_i,x\rangle^*\langle e_i,x\rangle$
converges in norm. In particular, the Fourier coefficients
$\{{\langle e_i,x\rangle}\}_{i\in I}$ of $x$ belong to $H_{A,I}$.
Thus, one has a well defined $A$-linear map from $M$ to $H_{A,I}$
given by the rule:
\[
x\mapsto \{{\langle e_i,x\rangle}\}_{i\in I}.
\]
The straightforward verification shows that this map is actually
an isomorphism. This result was extended for the frame context in
\cite[Theorem 4.1]{FrankLarJOP}.

Moreover the cardinality of an orthonormal basis in a Hilbert
module is unique like it happens for a Hilbert space. Indeed, as
it happens for a Hilbert space, the norm convergence of the series
$\sum\nolimits_{i\in I}a_i^*a_i$ implies that the number of its
non-zero entries are at most countable. And then it remains only
to apply well-known arguments similar to the ones for the Hilbert
space case (cf. \cite[I.\S 5.4]{Naimark_NR}).

\begin{prop} Any two closed orthonormal systems of a Hilbert module
over a unital $C^*$-algebra have the same cardinality. $\;\Box$
\end{prop}

Let us remark that the cardinality of a complete quasi-orthonormal
system  in a von Neumann module is not unique \cite[Remark
4.15]{SkeideRIA}. It is easy to see that the same is true for
closed quasi-orthonormal systems (for instance, we can consider
functions $f_1, f_2$ and $1$ of Example
\ref{ex:orthogonal_basis}). 

The next example shows that there are orthonormal systems in
standard Hilbert $C^*$-modules that cannot be extended to complete
orthonormal systems, a situation that differs from the cases of
orthogonal systems described in Proposition
\ref{teo:full_orth_syst}.  A natural question is the existence of
examples for a separable algebra $A$.

\begin{ex}\rm Assume $A=L^\infty [0,1]$, $M=l_2(A)$ and choose the functions
$f_1$, $f_2$ as in Example \ref{ex:orthogonal_basis}. Let
$\{e_i\}_{i=1}^\infty$ be the standard basis of $l_2(A)$
meaning that all entries of $e_i$ are zero except the $i$-th, which is
the identity of $A$. Then the vectors $\{x_i\}_{i=1}^\infty$ of
$M$, where $x_i=f_1e_i+f_2e_{i+1}$, form an
orthonormal system. It is not complete; indeed, suppose a
vector $y=(g_1,g_2,\dots)$ of $M$ is orthogonal to $x_i$ for any
$i$, that is its entries are such that:
\begin{gather}\label{eq:ex_ONS_l2A}
g_1|_{[0,1/2]}=0, \quad g_i=0 \,\,\text{for}\,\, i\ge 2 ;
\end{gather}
this holds, for instance, for the non-zero vector $x=f_2e_1$. On
the other hand the family $\{x_i\}_{i=1}^\infty$ cannot be
enlarged to a complete orthonormal system, because the inner
square of any vector  satisfying (\ref{eq:ex_ONS_l2A})  cannot
give the identity. Let us remark, by the way, that the vector $x$
extends the set $\{x_i\}_{i=1}^\infty$ to a complete orthogonal
system.
\end{ex}

The next example shows that there are complete orthonormal systems
in standard Hilbert $C^*$-modules which are not closed. This is
one of the crucial differences between general Hilbert
$C^*$-modules and Hilbert spaces.

\begin{ex}\label{ex:ONS_full_but_not_closed}\rm (This example was refined with crucial suggestions from M. Skeide). Suppose $A=L^\infty[0,1]$ and $M=l_2(A)$ is the standard
countably generated module over $A$. The desired system
$\{e_i\}_{i=1}^\infty$ of $M$, where
$e_i=(f_{i1},f_{i2},f_{i3},\dots)$ is constructed as follows.
 Let us denote by $\varphi_{[a,b]}$ the characteristic
function of the interval $[a,b]$, i.e.
\[
\varphi_{[a,b]}(x)=
\begin{cases}
1, & \text{if $\,\, x\in [a,b]$;}\\
0, & \text{otherwise}.
\end{cases}
\]
Consider $c_i=1-\frac{1}{2^i}$ for any non-negative integer $i$.
Then
\begin{gather*}
f_{i1}=\varphi_{[c_{i-1},c_i]},\quad i\ge 1,\\
f_{i(i+1)}=\varphi_{[c_i,1]},\quad i\ge 1,\\
f_{ii}=\varphi_{[0,c_{i-1}]},\quad i>1,
\end{gather*}
and $f_{ij}=0$ for all other positive integer values of $i$ and
$j$.

For such a construction we have the following
properties:
\begin{enumerate}
    \item only a finite number of functions $\{f_{ij}\}_{j=1}^\infty$ is
    non-zero for any $i$, apart from this, the sum $\sum_{j=1}^\infty f_{ij}=1$
    everywhere on the interval $[0,1]$
    (except either the points $c_{i-1}$ and $c_i$ if $i\ge 2$ or the
    point $c_1$ if $i=1$, but subsets of
zero measure are not significant). This implies: $\langle
e_i,e_i\rangle=1$ for any $i$;

    \item whenever $i\neq k$ the supports of the functions
    $f_{ij}$ and $f_{kj}$ do not intersect each other for any $j$.
    This implies: $\langle e_i,e_k\rangle=0$ for  $i\neq k$;

    \item for any $i$ the union over $j$ of the supports of the functions
    $f_{ij}$ coincides with the
    interval $[0,1]$;  this means that the system $\{e_i\}$ is complete.
\end{enumerate}
But at the same time the system $\{e_i\}$ cannot be closed
since, for example, for the vector $x=(1,0,0,\dots)$ the series
$\sum\nolimits_{i=1}^\infty \langle x,e_i\rangle\langle e_i,x\rangle$ does
not converge in norm.
\end{ex}

We remark that a situation similar to the one of the previous
example cannot occurs for finite orthonormal systems since,
clearly, any finite complete orthonormal system of a Hilbert
$C^*$-module is closed. Now we would like to describe an example
of a countably, but not finitely generated Hilbert $C^*$-module
possessing a finite complete orthogonal system. In fact the idea of
the next example may be used for constructing families of
such modules, corresponding to the branched coverings over compact
Hausdorff spaces.

\begin{figure}[h]
\begin{picture}(90,60)
\put(0,0){\line(1,0){70}} \put(80,5){$X$}
\put(0,40){\line(1,0){40}} \put(40,40){\line(3,2){30}}
 \put(40,40){\line(3,-2){30}}  \put(80,50){$Y$}
\put(30,30){\vector(0,-1){15}} \put(20,20){$p$}
\end{picture}
 \caption{Example~\ref{ex:count_gen_mod}}\label{fig:count_gen_mod}
\end{figure}

\begin{ex}\label{ex:count_gen_mod}\rm
Let us consider the map $p: Y\rightarrow X$ from Figure
\ref{fig:count_gen_mod}, where $X$ is an interval, say $[-1,1]$,
and $Y$ is the topological union of one interval with two copies
of another half-interval with a branch point at $0$. Then $C(Y)$
is a Banach $C(X)$-module for the action:
\[
(f\xi)(y)=f(y) \, \xi(p(y)), \qquad \mathrm{for} \quad f\in C(Y), \quad \xi\in C(X).
\]
Let us define the $C(X)$-valued inner product on $C(Y)$ by the
formula
\begin{gather}\label{eq:inner_pr}
    \langle f,g\rangle (x)=\frac{1}{\# p^{-1}(x)} \sum\nolimits_{y\in
    p^{-1}(x)} \overline{f(y)}g(y),
\end{gather}
where $\# p^{-1}(x)$ is the cardinality of $p^{-1}(x)$. It was
shown in  \cite{PavTro_cov} that $C(Y)$ is a countably, but not
finitely generated Hilbert $C(X)$-module with respect to the inner
product (\ref{eq:inner_pr}). The space $Y$ consists of the three
intervals with the common boundary point; interval that we number
in some arbitrary way. Then, for $i=1, 2, 3$, let us consider all
continuous on $Y$ functions $f_i$ that are not zero at all points
of the $i$-th interval except the boundary point and are zero at
the others points of $Y$. Clearly, the functions $f_1, f_2, f_3$
form a finite orthogonal complete system of $C(Y)$.
\end{ex}

We finish the paragraph by describing one family of non-standard
bases for $l_2(L^\infty[0,1])$.

\begin{ex}\rm
Let $A=L^\infty[0,1]$ and $M=l_2(A)$. For any
positive integer $n$, consider the functions
\[
f_i(x)=
\begin{cases}
1, & \text{if $\,\, \frac{i-1}{n}\le x < \frac{i}{n}$;}\\
0,& \text{otherwise}
\end{cases}
\]
and the matrix
\[F_n=
\begin{pmatrix}
f_1 & f_2 & \dots & f_{n-1}&f_n \\
f_2 & f_3 & \dots & f_n&f_1\\
 \dots &  \dots & \dots &  \dots &  \dots \\
f_{n} & f_1 & \dots & f_{n-2}& f_{n-1}
\end{pmatrix}
\]
We form a new matrix with an infinite number of rows
and columns in the following manner:
\[B=
\begin{pmatrix}
F_n & 0&\dots & 0&\dots \\
0 & F_n&\dots & 0&\dots\\
 \dots &  \dots & \dots &  \dots &  \dots\\
0 & 0&\dots & F_n& \dots \\
 \dots &  \dots & \dots &  \dots &  \dots
\end{pmatrix} ,
\]
and introduce vectors $e_i$ as the i-th rows of $B$.
The system $\{e_i\}$ forms an orthonormal basis of the Hilbert module $l_2(A)$.
\end{ex}

{\bf Acknowledgement:}
We are grateful to M. Frank and M. Skeide for helpful
remarks and discussions.
The work was partially supported by the `Italian project Cofin06 - Noncommutative geometry,
quantum groups and applications'. AP partially supported by RFBR (grant 08-01-00034) and
by the joint RFBR-DFG project (grant 07-01-91555). This work was completed  at  the Chern Institute of Mathematics, Nankai University, Tianjin, China; we are grateful for the nice hospitality there.

\bibliographystyle{plain}

\end{document}